\documentclass[a4,12pt,reqno,fleqn]{amsart}
\usepackage{graphicx}
\usepackage{amscd}
\usepackage{xfrac}
\usepackage{mathtools}
\emergencystretch=\maxdimen
\allowdisplaybreaks[1]
\usepackage{array}
\usepackage{xcolor}
\usepackage[all]{xy}
\subjclass[2020]{46B20, 46E40, 46G10}

\keywords{Birkhoff-James orthogonality, symmetric points, smooth points, Fr$\acute{e}$chet derivative, Bochner integrable functions}

\usepackage{amsmath}
\usepackage{amsthm}
\usepackage{amsfonts,amssymb,hyperref,mathrsfs,cleveref,setspace,enumitem,verbatim}
\usepackage[margin=3.5cm]{geometry}

\DeclareMathAlphabet{\mathpzc}{OT1}{pzc}{m}{it}

\newtheorem{thm}{Theorem}[section]
\newtheorem{cor}[thm]{Corollary}
\newtheorem{lem}[thm]{Lemma}

\theoremstyle{definition}

\newtheorem{remark}[thm]{Remark}
\newtheorem{example}[thm]{Example}

\newcommand{\bj}{\perp_{BJ}}
\newcommand{\nbj}{\not\perp_{BJ}}

\author[Mohit and R.~Jain]{Mohit and Ranjana Jain}

\address{Mohit, Department of  Mathematics, University of Delhi, Delhi}
\email{mohitdhandamaths@gmail.com}

\address{Ranjana Jain, Department of Mathematics, University of Delhi, Delhi}
\email{rjain@maths.du.ac.in}
\thanks{Research of the first named author is supported by Savitribai Jyotirao Phule Single Girl Child Fellowship vide F.No. 82-7/2022(SA-III)}
\begin{document}
		\title{Some geometric properties of spaces of vector-valued integrable functions}
	\maketitle
	\textbf{Abstract:} We identify the smooth points of $L^1(\mu,X)$, and provide some necessary and sufficient conditions for left and right symmetry of points with respect to Birkhoff-James orthogonality in $L^p(\mu,X), 1\leq p<\infty$, where $\mu$ is any complete positive measure  and $X$ is a Banach space with some suitable properties. 
		\section{Introduction}
	Birkhoff-James orthogonality (in short B-J orthogonality) plays an important role in the geometry of Banach spaces. An element $x$ in a normed space $X$ over $\mathbb{K}$ is said to be B-J orthogonal to  $y \in X$ (written as $x\bj y$) if
	$$\|x+\lambda y\|\geq\|x\|,\; \text{for all} \;\lambda\in \mathbb{K}.$$ 
	Unlike the usual orthogonality in the Hilbert spaces, B-J orthogonality is not symmetric. To gain a better insight of the elements which fail symmetry,  the concept of local symmetry with respect to the B-J orthogonality is studied through the left and right symmetric points. To recall, an element $x$ in a normed space $X$ is said to be a {\it left symmetric point} (respectively, {\it right symmetric point}) if $x\bj y$ implies $y\bj x$ (respectively,, $y\bj x$ implies $x\bj y)$ for all $y\in X$. Nowadays, the study of local symmetry  is an active area of research due to its   applications in the geometry of Banach spaces. 	Many researchers carried out the study of local symmetry in various spaces,  for example, see \cite{sain, ghosh,komuro, saito, vitor,  bose1}. 
	
	It is evident that in many well known spaces the collection of left and right symmetric points is sparse. For example, there is no non-zero left symmetric point in the sequence space $\ell^1$ whereas only the scalar multiples of $e_n,\;n\in\;\mathbb{N}$ are the right symmetric points in $\ell^1$ (see \cite[Theorem 3.6, Theorem 3.7]{bose1}). The same pattern can be seen in $B(H)$, where zero is the only left symmetric point and the right symmetric points are precisely the scalar multiples of isometries or coisometries (see \cite[Corollary 3.4, Corollary 4.5]{turnsek}). For the commutative unital $C^*$-algebra $C(K)$, the space of continuous functions on a compact Hausdorff space $K$, the left symmetric points are precisely the functions which vanish everywhere except possibly at an isolated point, whereas the right symmetric points are the functions which attain their norms at every point (see \cite[Theorem 3.6, Theorem 4.6]{komuro}).  On the other hand, Deeb and Khalil, in 1994, characterized smooth points in $L^1(\mu,X)$ for a finite measure $\mu$ and a Banach space $X$ with separable dual  (see, \cite[Theorem 1.1]{deeb}). Further, Bose \cite[Theorem 3.0.3]{bose} proved that a point $f\in L^1(\mu)$ is smooth if and only if $f$ is non-zero almost everywhere.

	The aim of this article is to analyze similar questions in a more general setting of the space of vector-valued functions. An important tool for the discussion of local symmetry in the spaces $L^p(\mu,X),\;1\leq p<\infty$ will be the characterizations of the B-J orthogonality in these spaces derived by the authors recently (see \cite[Corollary 3.5, Corollary 3.7]{jain}). 
	
	
	Here is a quick overview of the main results of this article. Throughout the article, we assume $(S,\mu)$ to be a complete positive measure space.
	In Section 3, we first prove that for a Banach space $X$ whose norm is  Fr$\acute{e}$chet differentiable, the smooth points in $L^1(\mu,X)$ are  precisely the functions which are non-zero almost everywhere.  We then establish that approximately smooth and smooth points coincide in $L^1(\mu,X)$. Next, for a $\sigma$-finite measure $\mu$ and a complex smooth Banach space or a real Banach space whose norm is Fr$\acute{e}$chet differentiable, we prove that if a non-zero element $f\in L^1(\mu,X)$ is a left symmetric point, then either $Z(f)^c$ is an atom or is a disjoint union of two atoms, where $Z(f):=\{s\in S: f(s) =0\}$. We provide examples to show that the converse of this statement fails. A sufficient condition for $f\in L^1(\mu,X)$ to be a left symmetric point is also given. We then provide some necessary and sufficient conditions for right symmetry of elements in  $L^1(\mu,X)$ in terms of atoms. We conclude this section by proving similar results of local symmetry in the space $L^p(\mu,X),\;1<p<\infty,\;p\neq2$. These conditions lead to the fact that for a non-atomic measure $\mu$, the space  $L^p(\mu,X), 1\leq p<\infty,\;p\neq2$ has no non-zero left and right symmetric points. 
	
	\section{Preliminaries}
	
	We first collect some definitions and results which we need for our investigation.
	For a measure space $(S,\mu)$ and a Banach space $X$, we shall consider the spaces  
	$$	 L^{p}(\mu, X)= \bigg\{f:S\rightarrow X|\ f \ \text{is strongly measurable and} \;\int\limits\limits\limits_{S}\|f(s)\|^p\,ds <\infty \bigg\},$$ 
	when	$1\leq p<\infty,$ known as the Lebesgue-Bochner space;
	and 
	$$ L^{\infty}(\mu, X)= \{f:S\rightarrow X|\ f \ \text{is strongly measurable and}  \;\text{ess sup}\;\|f(s)\| <\infty \},$$
	where almost everywhere equal functions are identified in both the spaces.
	
	For a non-zero element $ x\in X$, a {\it support map}  at $x$ is a bounded linear functional $F_{x}$ on $X$ of norm one satisfying $F_{x}(x)=\|x\|.$  An element $x$ is said to be {\it smooth} if $F_x$ is unique. The space $X$ is {\it smooth} if every non-zero element of $X$ is smooth. It is known that $x \bj y$ if and only if  there exists a support map $F_x$ such that $F_x(y) =0$.
	
	A norm function $\| \cdot \|$ on $X$ is said to be {\it Fr$\acute{e}$chet differentiable} at a non-zero element $ x\in X$ if there exists an $f \in X^*$  satisfying $$\lim_{h \to 0}\frac{\big|\|x+h\|-\|x\|-f(h)\big|}{\|h\|}=0.$$ The norm function is Fr$\acute{e}$chet differentiable if it is Fr$\acute{e}$chet differentiable at every non-zero point. It is known that $X$ is smooth if its norm is Fr$\acute{e}$chet differentiable. We will use the following characterizations quite often: 	
	\begin{thm}\cite [Corollary 3.5]{jain}\label{ourbj}
		Let $X$ be a Banach space and $f,g\in L^{1}(\mu,X)$. Then $f\perp_{BJ}g$ in  $L^{1}(\mu,X)$ if and only if 
		\begin{equation}\label{eqa}
			\hspace{1.5cm}		\Bigg|\int\limits\limits\limits_{Z(f)^c}{Re(F_{f(s)}( \alpha g(s))) \,ds}\Bigg|\leq\int\limits\limits\limits_{Z(f)} |\alpha| \| g(s)\| \,ds,\ \forall \ \alpha\in\mathbb{C},
		\end{equation} 
		when $X$ is a complex smooth Banach space, or 
		\begin{equation}\label{eqb}
			\hspace{3cm}	\Bigg|\int\limits_{Z(f)^c}{F_{f(s)}(g(s)) \,ds}\Bigg|\leq\int\limits\limits\limits_{Z(f)} \| g(s)\| \,ds,
		\end{equation}
		when $X$ is a real Banach space whose norm is Fr$\acute{e}$chet differentiable. 	
	\end{thm}
	\begin{thm}\cite[Corollary 3.7]{jain}\label{ourbjp}
		Let $X$ be a Banach space whose norm is Fr$\grave{e}$chet differentiable. For $f,g\in L^{p}(\mu,X)$, $1<p<\infty$,  $f\perp_{BJ}g$ in  $L^{p}(\mu,X)$ if and only if 		
		$$\int\limits_{Z(f)^c}{\|f(s)\|^{p-1}F_{f(s)}(g(s))\,ds}=0.$$
	\end{thm}
\section{Geometric properties of Lebesgue-Bochner space}
We will first determine the smooth and approximate smooth points of $L^{1}(\mu,X)$. Note that the class of Banach spaces with separable dual (as discussed by Deeb and Khalil) is different from the class of Banach spaces whose norm is Fr$\acute{e}$chet differentiable.  

\begin{thm}\label{smooth point} 
	Let $X$ be a Banach space whose norm is Fr$\acute{e}$chet differentiable. An element $f\in L^{1}(\mu,X)$ is smooth if $f\neq0$ almost everywhere. The converse is true if, in addition, $\mu$  is a $\sigma$-finite measure. 
\end{thm}
\begin{proof}
	First, assume that $f\neq0$ almost everywhere, so that $\mu(Z(f))=0$. It is sufficient to show that  B-J orthogonality is right additive at $f$ (see \cite[Theorem 2.3.2, Remark 2.3.4]{mal}). Consider $g_1,\; g_2\in\;L^{1}(\mu,X)$ such that $f\bj g_{1}$ and $f\bj g_{2}$. If $X$ is a complex Banach space, then by \Cref{ourbj}, for any $\alpha\in\mathbb{C}$, 
	$$		\int\limits\limits\limits_{Z(f)^c}{Re(F_{f(s)}( \alpha g_i(s))) \,ds}=0, \ \ \text{for} \ i=1,2.
	$$  
	Using linearity of support map we have
	\begin{equation*}
		\hspace{1.5cm}		\int\limits\limits\limits_{Z(f)^c}{Re(F_{f(s)}( \alpha (g_1+g_2)(s))) \,ds}=0,
	\end{equation*}
	which gives $f\bj (g_1+g_2)$, so that right additivity holds at $f$. The case when $X$ is a real Banach space goes on the similar lines.
	
	For the converse, let if possible, $\mu(Z(f))>0$. Since $f$ is non-zero, there exists $s_0\in S$ such that $f(s_0)\neq0$. Define  $h_{i}:S\rightarrow X^{*}$, $i=1, 2$ as 
	$$
	h_1(s)= \Bigg\{
	\begin{array}{l l}
		F_{f(s)} & \text{if}\;s\in  Z(f)^c\\
		0 & \text{if}\;s\in Z(f)
	\end{array}\Bigg\}, \,\,
	\quad
	h_2(s)=\Bigg\{
	\begin{array}{l l}
		F_{f(s)} & \text{if}\;s\in  Z(f)^c\\
		F_{f(s_0)}	 & \text{if}\;s\in Z(f)
	\end{array}\Bigg\}.
	$$
	We claim that $h_{i}\in L^{\infty}(\mu,X^{*})$ for $i=1,2$. For this, let $\{f_n\}$ be a sequence of simple measurable functions that converges to $f$ almost everywhere. Define the sequences $\{\phi_n\}$ and $\{\psi_n\}$ as 
	$$
	\phi_n(s)= \Big\{
	\begin{array}{l l}
		F_{f_n(s)} & \text{if}\; s\in  Z(f)^c\cap Z(f_n)^c\\
		0 & \text{otherwise}
	\end{array}\Big\}, $$
	$$\psi_n(s)= \Big\{
	\begin{array}{l l}
		F_{f_{n}(s)} & \text{if}\; s\in  Z(f)^c\cap Z(f_n)^c\\
		F_{f(s_0)} &  \text{otherwise}
	\end{array}\Big\}.
	$$
	Since the norm on $X$ is Fr$\acute{e}$chet differentiable, the map $x \mapsto F_x$ is continuous on $X \setminus \{0\}$ (see, \cite[Proposition 3.2]{jain}), so it is easy to check that both $\{\phi_n\}$ and $\{\psi_n\}$ are sequences of simple measurable functions converging to $h_{1}$ and $h_{2}$ almost everywhere, respectively. Thus, $h_1$ and $h_2$ are measurable, the measure $\mu$ being complete. Clearly, $h_{i}\in L^{\infty}(\mu,X^{*})$ with $\|h_i\|_{\infty}=1$, for $i=1,2$.
	Consider the maps $T_{h_i}:L^{1}(\mu,X)\rightarrow \mathbb{C}$ given by $T_{h_i}(f)=\int\limits_{S}h_i(s)(f(s))\,ds$ for $i=1,2$. Now, $\mu(Z(f))>0$ and $\mu$ is $\sigma$-finite, therefore there exists $A\subseteq Z(f)$ such that $0<\mu(A)<\infty$.  Consider $p \in L^1(\mu,X)$  defined as $p=f(s_0)\chi_{A}$. Then, $T_{h_1}(p)=0$ and $T_{h_2}(p)=\|f(s_0)\|\mu(A)\neq 0$, so that $T_{h_i}$'s are different support maps at $f$, which contradicts the fact that $f$ is smooth. Hence $\mu(Z(f))=0$, which gives $f\neq0$ almost everywhere.    
\end{proof}
\begin{remark}
	It is worth to mention that if $X$ is smooth real Banach space, then  $L^p(\mu, X), 1 < p < \infty$ is smooth (see \cite[Theorem 3.1]{leo}). Also, a complex normed space is smooth if and only if it is smooth as a real normed space. Thus, for a smooth Banach space $X$, all the non-zero elements of $L^p(\mu, X), 1 < p < \infty$ are smooth.
\end{remark}
Recently Chmieli$\acute{n}$ski et. al. \cite{jacek}, introduced the concept of approximate smoothness of elements in Banach spaces. For a normed space $X$, a non-zero element $ x\in X$ is said to be {\it approximately smooth} if diam$J(x) < 2$,  where $J(x)$ is the set of all support functionals at $x$  and diam$J(x)=\sup \{\|F_x-G_x\|: F_x,\; G_x\in J(x)\}$. Clearly every smooth point is approximately smooth, but the converse need not be true in general. It is interesting to note that approximate smoothness coincides with smoothness in $L^1(\mu,X)$ for an appropriate class of Banach spaces and measure spaces. 
\begin{cor}
	Let $X$ be a Banach space whose norm is Fr$\acute{e}$chet differentiable and  $\mu$ be a $\sigma$-finite measure. Then $f\in L^1(\mu,X)$ is approximately smooth if and only if  $f\neq0$ a.e.
\end{cor}
\begin{proof}
	Let $f\in L^1(\mu,X)$ be an approximately smooth point. Let if possible, $\mu(Z(f))\neq0$. Let $s_0\in S$ be a point for which $f(s_0)\neq0$. Define  $g_{i}:S\rightarrow X^{*}$, $i=1, 2$ as 
	$$
	g_1(s)= \Bigg\{
	\begin{array}{l l}
		F_{f(s)} & \text{if}\;s\in  Z(f)^c\\
		F_{f(s_0)} & \text{if}\;s\in Z(f)
	\end{array}\Bigg\}, \,\,
	\quad
	g_2(s)=\Bigg\{
	\begin{array}{l l}
		F_{f(s)} & \text{if}\;s\in  Z(f)^c\\
		-F_{f(s_0)}	 & \text{if}\;s\in Z(f)
	\end{array}\Bigg\}.
	$$
	As done previously, since the norm on $X$ is Fr$\acute{e}$chet differentiable, clearly $g_i\in L^{\infty}(\mu, X^{*})$ with $\|g_i\|_{\infty}=1$. Define $T_{g_i} \in (L^1(\mu, X))^*$ as $T_{g_i}(h)=\int\limits_{S}g_{i}(s)(h(s))\,ds$, then $T_{g_i} \in J(f)$ for $i=1,2$. Since the measure $\mu$ is $\sigma$-finite and $\mu(Z(f))>0$, therefore there exists $A\subseteq Z(f)$ such that $0<\mu(A)<\infty$. Let $h:S\rightarrow X$ be defined as $ h=f(s_0)\chi_{A}$. Clearly, $h\in L^1(\mu,X)$ and  $T_{g_1}(h)=\|f(s_0)\|\mu(A)=\|h\|=-T_{g_2}(h)$. Thus, 
	$$\|T_{g_1}-T_{g_2}\|\geq \frac{|T_{g_1}(h)-T_{g_2}(h)|}{\|h\|} =2,$$
	which gives that diam$J(f)\geq2$. Hence, $f$ is not approximately smooth.
	
	Converse follows directly from \Cref{smooth point}.
\end{proof}
We next discuss the left symmetric points of the Bochner space. We would like to mention that our results related to local symmetry in $L^p(\mu,X),1\leq p<\infty,\;p\neq2$ are motivated from those in $L^p(\mu), 1\leq p<\infty,\;p\neq2$, the space of scalar-valued integrable functions,  given by Bose \cite{bose}. 
Let us recall that for a measure space $(S,\mu)$,  a measurable set $A$ is called an {\it atom} if $\mu(A)>0$ and for any measurable subset $B$ of $A$, either $\mu(B)=0$ or $\mu(A)=\mu(B)$ (see, \cite{dunford}).  Note that if an atom $A$ is $\sigma$-finite then $\mu(A)<\infty$ so that for any  measurable subset $B$ of $A$, either $\mu(B)=0$ or $\mu(A\setminus B)=0$. The space is {\it non-atomic} if there are no atoms.

\begin{thm}\label{left symmetry}
	Let $X$ be either a complex smooth Banach space or a real Banach space whose norm is Fr$\acute{e}$chet differentiable, and $\mu$ be a $\sigma$-finite measure. If $ f\in L^1(\mu,X)$ is a non-zero left symmetric point, then  $Z(f)^c$ is either an atom or is a disjoint union of exactly two atoms.
\end{thm}
\begin{proof}
	We prove the result for the complex smooth Banach space $X$ and the proof for the real Banach space whose norm is Fr$\acute{e}$chet differentiable goes on the similar lines using (\ref{eqb}). We first claim that $\mu(Z(f))=0$. If not, then there exists a measurable subset, say $A$ of $Z(f)$, such that $0<\mu(A)<\infty$, since $\mu$ is $\sigma$-finite. For a fixed unit vector $x_0\in X$, define   $h=f \chi_{Z(f)^c}+ \frac{\|f\|}{\mu(A)}x_0 \chi_A$. Clearly, 
	$h\in L^1(\mu,X)$ and for any scalar $\alpha\in \mathbb{C}$,
	$$\Bigg|\int\limits_{Z(f)^c}Re(F_{f(s)}(\alpha h(s)))\,ds\Bigg| = |Re(\alpha)|\|f\|\leq \int\limits_{Z(f)}|\alpha|\|h(s)\|\,ds, $$ 
	so that, by \Cref{ourbj}, $f\bj h$. However, for $\alpha\in \mathbb{C}$ for which Re$(\alpha)\neq0$ we have 
	\begin{align*}
		\hspace*{-1.3cm} 	\Bigg|\int\limits_{Z(h)^c}Re(F_{h(s)}(\alpha f(s)))\,ds\Bigg|& = \Bigg|\int\limits_{Z(f)^c}Re(F_{f(s)}(\alpha f(s)))\,ds + \int\limits_{A}Re(F_{h(s)}(\alpha f(s)))\,ds \Bigg|\\
		&= |Re(\alpha)|\|f\|> \int\limits_{Z(h)}|\alpha|\|f(s)\|\,ds=0.
	\end{align*}
	Thus,  $h\nbj f$, which contradicts the fact that $f$ is left symmetric and hence the claim. \\
	Now, let us consider that neither $Z(f)^c$ is an atom nor $Z(f)^c$ is a disjoint union of exactly two atoms. Then, by the definition of atom, one can find three disjoint measurable subsets say $P,\;Q$ and $R$ of $Z(f)^c$ with positive measures. Now, as $f$ is non-zero on $P,\;Q$ and $R$, therefore either
	
	\begin{equation}\label{equation1}
		0<\int\limits_{Q\cup R}\|f(s)\|\,ds\leq\int\limits_{P}\|f(s)\|\,ds
	\end{equation}
	or
	\begin{equation}\label{eq2} 
		\int\limits_{Q\cup R}\|f(s)\|\,ds>\int\limits_{P}\|f(s)\|\,ds>0.
	\end{equation}
	If (\ref{equation1}) is true, then using the fact that $P\cup R\subseteq Q^c$, we further get  $$0<\int\limits_{Q}\|f(s)\|\,ds<\int\limits_{P\cup R}\|f(s)\|\leq\int\limits_{Q^c}\|f(s)\|\,ds.$$ 
	Consider $\beta_1=\int\limits_{Q^c}\|f(s)\|\,ds = \int\limits_{Z(f)^c \setminus Q}\|f(s)\|\,ds $ and $\beta_2=\int\limits_{Q}\|f(s)\|\,ds$, then $\beta_1 > \beta_2 >0$. Now, define $g\in\;L^1(\mu,X)$ as $g= \beta_1 f \chi_Q - \beta_2 f \chi_{Q^c}$.
	Then
	\begin{align*}
		\hspace*{-1.5cm}	\Bigg|\int\limits_{Z(f)^c}Re(F_{f(s)}(\alpha g(s)))\,ds\Bigg|&=|Re(\alpha)|\Bigg|\beta_1\int\limits_{Q\cap Z(f)^c}\|f(s)\|\,ds-\beta_2\int\limits_{Q^c\cap Z(f)^c}\|f(s)\|\,ds\Bigg|\\
		&=|Re(\alpha)(\beta_1\beta_2-\beta_2\beta_1)|\\&= 0 =\int\limits_{Z(f)}\|\alpha g(s)\|\,ds,
	\end{align*}
	for all $\alpha\in \mathbb{C}$, which implies that $f \bj g$ using \Cref{ourbj}.  Further, since $X$ is smooth, $F_{\alpha x}=F_{x}$ if $\alpha>0$ and $F_{\alpha x}=-F_x$ if $\alpha<0$. Also $Z(f)^c=Z(g)^c$, thus for $\alpha\in\mathbb{C}$ for which Re$(\alpha)\neq0$,
	\begin{align*}
		\hspace*{-0.75cm}	\Bigg|\int\limits_{Z(g)^c}Re(F_{g(s)}(\alpha f(s)))\,ds\Bigg|& =|Re(\alpha)|\Bigg|\int\limits_{Q\cap Z(g)^c}\|f(s)\|\,ds-\int\limits_{Q^c\cap Z(g)^c}\|f(s)\|\,ds\Bigg|\\
		&=|Re(\alpha)|(\beta_1-\beta_2)\\
		&	> 0=\int\limits_{Z(g)}\|f(s)\|\,ds.
	\end{align*}
	so that $g \nbj f$, which is a contradiction. However, if \ref{eq2} holds, then we have $0<\int\limits_{P}\|f(s)\|\,ds<\int\limits_{P^c}\|f(s)\|\,ds$ and we finish as before with $Q$ replaced by $P$.
\end{proof}
\begin{cor}\label{lspsmooth}
	For a $\sigma$-finite measure $\mu$ and  a Banach space $X$  whose norm is Fr$\acute{e}$chet differentiable, every non-zero left symmetric point in $L^1(\mu,X)$ is smooth.
\end{cor}
\begin{proof}
	For a non-zero $f \in L^1(\mu,X)$, as done in the proof of \Cref{left symmetry}, $\mu(Z(f))=0$. Thus, by \Cref{smooth point}, $f$ is a smooth point. 
\end{proof}
\begin{cor}
	For a Banach space $X$ whose norm is Fr$\acute{e}$chet differentiable, the sequence space $\ell^1(X)$ has no nonzero left symmetric point.
\end{cor}

In general, the converse of \Cref{left symmetry} is not true. Let us look at some examples to demonstrate this.
\begin{example}\label{examatom} $f$ need not be a left symmetric point if $Z(f)^c$ is an atom. 
	For this, consider the measure space $(S,\Sigma,\mu)$ defined as $S=(0,1)\cup(2,3)$, $\Sigma=\{\emptyset,\;S,\;(0,1),\;(2,3),\;A, \ (0,1)\cup A:A \subseteq (2,3)\}$, and $\mu$ is defined as $\mu(\emptyset)=0$ $=\mu((2,3))$ and $\mu((0,1)\cup A) =1$, for all $A\subseteq (2,3)$. Consider the infinite dimensional Banach space  $\ell^p,\;p>2$ whose norm is Fr$\acute{e}$chet differentiable. Since B-J orthogonality is not symmetric in $\ell^p$ \cite[Theorem 1]{jamesips},  there exist non-zero elements $x_0,\;y_0\in \ell^p$ such that $x_0\bj y_0$ but $y_0\nbj x_0$. Define $f,\;g \in L^1(\mu, \ell^p)$ as $f=x_{0}\chi_{(0,1)}$ and $g=y_{0}\chi_{(0,1)}$. Then
	$$\Bigg|\int\limits_{(0,1)}F_{f(s)}(g(s))\,ds\Bigg|=|F_{x_0}(y_0)|=0 = \int\limits_{(2,3)} \|g(s)\| ds,$$ which gives $f \bj g$ by \Cref{ourbj}.  However, $$\Bigg|\int\limits_{(0,1)}F_{g(s)}(f(s))\,ds\Bigg|=|F_{y_0}(x_0)|> 0 = \int\limits_{(2,3)} \|f(s)\| ds,$$ so that $g\nbj f$. Clearly, $Z(f)^c$ is an atom but $f$ is not a left symmetric point. 
\end{example}
\begin{example}\label{examtwoatom}
	If $Z(f)^c$ is a disjoint union of exactly two atoms, then $f$ need not be left symmetric. Let $S=(0,1)\cup(2,3)$ with $\Sigma=\{\emptyset,\;(0,1),\;(2,3),\;S\}$ and measure $\mu$ on $(S,\;\Sigma)$ be defined as $\mu(\emptyset)=0,\;\mu((0,1))=1=\mu((2,3)),\;\mu(S)=2$.
	For $	e_1,\;e_2 \in \ell^2$, define $f,\;g \in L^1(\mu, \ell^2)$ as $f=e_{1}\chi_{(0,1)}-e_{1}\chi_{(2,3)}$ and $g=(e_1+e_{2})\chi_{(0,1)}+e_{1}\chi_{(2,3)}$.	Clearly, $Z(f)^c$ is the disjoint union of exactly two atoms. Now,  by Riesz representation theorem, for $x \in \ell^2$, the support map $F_{x}$ on the Hilbert space $\ell^2$ is given by $F_{x}(\cdot)=\langle \cdot, \frac{x}{\|x\|}\rangle$. Therefore
	\begin{align*}
		\Bigg|\int\limits_{Z(f)^c}F_{f(s)}(g(s))\,ds\Bigg|&=|F_{e_1}(e_1+e_2)+F_{-e_1}(e_1)| \\
		& =\bigg|\langle e_1+e_2,\frac{e_1}{\|e_1\|}\rangle-\langle e_1,\frac{e_1}{\|e_1\|}\rangle\bigg|\\
		& = 0 = \int\limits_{Z(f)}\|g(s)\|\,ds,
	\end{align*}
	which using \Cref{ourbj} gives $f\bj g$. However,
	\begin{align*}
		\Bigg|\int\limits_{Z(g)^c}F_{g(s)}(f(s))\,ds\Bigg|&=|F_{e_1+e_2}(e_1)+F_{e_1}(-e_1)|\\ 
		& =\bigg|\langle e_1,\frac{e_1+e_2}{\|e_1+e_2\|}\rangle+\langle -e_1,\frac{e_1}{\|e_1\|}\rangle\bigg|\\ & =\bigg|\frac{1}{\sqrt{2}}-1\bigg|\\
		& > 0 = \int\limits_{Z(g)}\|f(s)\|\,ds,
	\end{align*}
	so that $g\nbj f$.  Hence $f$ is not a left symmetric point.
\end{example}

We next prove that a partial converse of \Cref{left symmetry} is true. In view of \Cref{lspsmooth}, it is evident that every non-zero left symmetric point in $L^1(\mu,X)$ is smooth. Thus, it is sufficient to check the converse only for smooth points.
\begin{thm}\label{leftsymL1}
	Let $X$ be a real Banach space whose norm is Fr$\acute{e}$chet differentiable and  $\mu$ be a $\sigma$-finite measure.	Let $f\in L^1(\mu,X)$ be a smooth point such that  $Z(f)^c$ is an atom and $f(s)$ is left symmetric for each $s\in S$, then $f$ is left symmetric.
\end{thm}
\begin{proof}
	Consider $g\in L^1(\mu,X)$ such that $f\bj g$. Since $f$ is smooth, by \Cref{lspsmooth}, $\mu(Z(f)) = 0$. Thus, from \Cref{ourbj},
	\begin{equation}\label{eq3} 
		\int\limits_{Z(f)^c}F_{f(s)}(g(s))\,ds= \int\limits_{Z(f)^c\cap Z(g)^c}F_{f(s)}(g(s))\,ds=0.
	\end{equation}
	Let $Z(f)^c=(Z(f)^c\cap Z(g)^c)\cup (Z(f)^c\cap Z(g)) = A \cup B$, say. Since $Z(f)^c$ is $\sigma$-finite and an atom, therefore, either $\mu(A)=0$ or $\mu(B)=0$. 
	
	If $\mu(A)=0$, then 
	$$\Bigg|\int\limits_{Z(g)^c}F_{g(s)}(f(s))\,ds\Bigg| = \Bigg|\int\limits_{A}F_{g(s)}(f(s))\,ds\Bigg| = 0 \leq \int\limits_{Z(g)}\|f(s)\|\,ds,$$
	so that  $g\bj f$, by \cite[Corollary 3.5]{jain}.  
	
	Now, if $\mu(B)=0$, write $A=P_1\cup P_2$, where $P_1=\{s\in A: F_{f(s)}(g(s))\geq0\}$ and $P_2=\{s\in A: F_{f(s)}(g(s))<0\}$. We first claim that the sets $P_i$'s are measurable. For this, let $\{f_n\}$ and $\{g_n\}$ be sequences of simple measurable functions converging to $f$ and $g$ almost everywhere, respectively. Define $h_n,\;h:S\rightarrow \mathbb{R}$ as 
	$$
	h_n(s)= \Big\{
	\begin{array}{l l}
		F_{f_n(s)}(g_{n}(s)) & \text{if}\; s\in Z(f_n)^c\cap Z(f)^c\\
		0 & \text{otherwise}
	\end{array}\Big\}, $$
	and
	$$h(s)= \Big\{
	\begin{array}{l l}
		F_{f(s)}(g(s)) & \text{if}\; s\in  Z(f)^c\\
		0 &  \text{otherwise}
	\end{array}\Big\}.
	$$
	Observe that $h_n(s)$ converges to $h(s)$ for a.e. $s\in S$. To see this, let $s\in S$ be such that $\lim f_n(s)=f(s)$ and $\lim g_n(s)=g(s)$. If $s\in Z(f)$, then clearly $h(s)=0=h_{n}(s)$ for all $n\in \mathbb{N}$. If $s\in Z(f)^c$, then by taking a tail of the sequence $\{h_n(s)\}$ (if required, to make $f_n(s) \neq 0$) and using the fact that  the map $x\mapsto F_x$ is continous on $X\setminus\{0\}$ (see, \cite[Proposition 3.2]{jain}), one can easily obtain that $\lim  h_n(s)=h(s)$. Thus,  $h$ is measurable being the a.e. limit of the sequence $\{h_n\}$ of simple measurable functions. Therefore, $P_1=h^{-1}[0,\infty)\cap A$ and $P_2=h^{-1}(-\infty,0)$ both are measurable.
	
	We next claim that $\mu(P_2)= 0$. If not,  then $\mu(P_1) = 0$ since $A$ is an atom being a  finite positive measure subset of an atom.   
	Also, by Equation (\ref{eq3}),
	$$0=\int\limits_{A}F_{f(s)}(g(s))\,ds=\int\limits_{P_2}F_{f(s)}(g(s))\,ds.$$ which contradicts the definition of $P_2$. Therefore  $\mu(P_2)=0$ and  by Equation (\ref{eq3}), $$0=\int\limits_{A}F_{f(s)}(g(s))\,ds=\int\limits_{P_1}F_{f(s)}(g(s))\,ds.$$ Thus, $F_{f(s)}(g(s))=0$ for a.e. $s\in P_1$ so that $f(s)\bj g(s)$ for a.e. $s\in P_1$. Since $f(s)$ is a left symmetric point for each $s\in S$, therefore $g(s)\bj f(s)$ for a.e. $s\in P_1$ that is $F_{g(s)}(f(s))=0$ for a.e. $s\in P_1$. Also, using the fact that  $\mu(P_2)=0=\mu(B)$ we have
	\begin{align*}
		\Bigg|\int\limits_{Z(g)^c}F_{g(s)}(f(s))\,ds\Bigg|&=\Bigg|\int\limits_{A}F_{g(s)}(f(s))\,ds\Bigg|\\
		& =\Bigg|\int\limits_{P_1}F_{g(s)}(f(s))\,ds\Bigg|\\
		&	=0\\
		&=\int\limits_{B}\|f(s)\|\,ds=\int\limits_{Z(g)}\|f(s)\|\,ds.
	\end{align*}
	Hence, by \Cref{ourbj}, $g\bj f$.
\end{proof}

We next discuss the right symmetric points of $ L^1(\mu,X)$.
\begin{thm}\label{rightsymL1}
	Let $X$ be a real Banach space whose norm is Fr$\acute{e}$chet differentiable and  $f\in L^1(\mu,X)$ be a non zero element. If $f$ is right symmetric, then $Z(f)^c$ is an atom. The converse is true, if in addition, $f(s)$ is a right symmetric point for each $s\in S$.
\end{thm}
\begin{proof}
	Let $f$ be a right symmetric point and let, if possible, $Z(f)^c$ be not an atom. So there exists a proper subset $B$ of $Z(f)^c$ such that $0<\mu(B)<\mu(Z(f)^c)$. Define $g \in L^1(\mu,X)$ as 
	$$	g=
	\begin{cases}
		f \chi_B & \text{if}\;\int\limits_{B}\|f(s)\|\,ds\leq\int\limits_{Z(f)^c\setminus B}\|f(s)\|\,ds\\
		f \chi_{ Z(f)^c\setminus B} & \text{if}\;\int\limits_{Z(f)^c\setminus B}\|f(s)\|\,ds < \int\limits_{B}\|f(s)\|\,ds.
	\end{cases}$$
	Clearly, in both cases  
	$$\Bigg|\int\limits_{Z(g)^c}F_{g(s)}(f(s))\,ds\Bigg|\leq\int\limits_{Z(g)}\|f(s)\|\,ds$$
	so that $g\bj f$, by \Cref{ourbj} which further gives $f\bj g$.
	Now, if $g = 	f \chi_B$, then    $$\Bigg|\int\limits_{Z(f)^c}F_{f(s)}(g(s))\,ds\Bigg|\leq\int\limits_{Z(f)}\|g(s)\|\,ds=0$$ that is, $$\int\limits_{B}\|f(s)\|\,ds=0$$ which cannot be true as $\mu(B)>0$ and $\|f(s)\|>0$ for all $s\in B$. Similarly, we get a contradiction for the other case and hence $Z(f)^c$ is an atom.
	
	Conversely, let $Z(f)^c$ be an atom and $f(s)$ be right symmetric for every $s\in S$. Consider a non-zero element $g\in L^1(\mu,X)$  such that $g\bj f$. Thus, by \Cref{ourbj}, we have	
	$$	 \Bigg|\int\limits_{Z(g)^c\cap Z(f)^c}F_{g(s)}(f(s))\,ds\Bigg|\leq\int\limits_{Z(g)\cap Z(f)^c}\|f(s)\|\,ds.	 $$
	Now, $Z(f)^c=A\cup B$, where $A=Z(f)^c\cap Z(g)^c$ and $B=Z(f)^c\cap Z(g)$.  Since $Z(f)^c$ is $\sigma$-finite (see, \cite[p.no. 191]{folland}) and  an atom so either $\mu(A)=0$ or $\mu(B)=0$.	 If $\mu(A)=0$ then $$0=\Bigg|\int\limits_{Z(g)^c\cap Z(f)^c}F_{f(s)}(g(s))\,ds\Bigg|\leq\int\limits_{Z(f)\cap Z(g)^c}\|g(s)\|,\,ds$$	so that $f\bj g$.	If $\mu(B)=0$, the proof goes on the similar lines of \Cref{leftsymL1}. 
\end{proof}

\begin{remark}
	In general,  if $Z(f)^c$ is an atom, then $f$ need not be right symmetric. One can easily verify the same with the help of \Cref{examatom}, that $g$ is not a right symmetric point but $Z(g)^c$ is an atom.
\end{remark}
\begin{cor}\label{l1rsp}
	Let $X$ be a real Banach space whose norm is Fr$\acute{e}$chet differentiable. Then, the right symmetric points of $\ell^1(X)$ are precisely of the form  $(0,..,0,x_0,0,..)$, where $x_0$ is a right symmetric point of $X$.
\end{cor}
\begin{proof}
	In view of \Cref{rightsymL1}, it is sufficient to check that if $x=(0,..,0,x_0,0,..)$ is a right symmetric point, then $x_0$ is also  right symmteric in $X$. Let $y_0\in X$ be such that $y_0\bj x_0$. Clearly, for $y=(0,..,0,y_0,0,...) \in \ell^1(X)$, using the definition of B-J orthogonality, $y \bj x$, so that $x \bj y$, which in turn gives $x_0\bj y_0$ and hence the result.	
\end{proof}
Next, we discuss local symmetry in $ L^p(\mu,X),\;1<p<\infty,\;p\neq2$. We would like to mention that the statements and the proofs are similar to the case when $p=1$ with some appropriate modifications. However, we present the proofs here for the sake of completeness. 

\begin{thm}\label{symmLp}
	Let $X$ be a Banach space whose norm is Fr$\acute{e}$chet differentiable. If $f\in L^p(\mu,X),\;1<p<\infty,\;p\neq2$ is a non-zero left symmetric or right symmetric point, then $Z(f)^c$ is either an atom or a disjoint union of exactly two atoms.
\end{thm}
\begin{proof}
	Let if possible, $Z(f)^c$ is neither an atom nor a disjoint union of two atoms.  Then, as done earlier in \Cref{left symmetry}, there are disjoint measurable subsets say $P,\;Q$ and $R$ of $Z(f)^c$ with positive measure such that  $$0<\int\limits_{Q}\|f(s)\|^p\,ds<\int\limits_{P\cup R}\|f(s)\|^p\,ds\leq\int\limits_{Q^c}\|f(s)\|^p\,ds.$$ 
	We claim that $f$ is neither left symmetric nor right symmetric. For this, consider $\beta_1=\int\limits_{Q^c}\|f(s)\|^p\,ds = \int\limits_{Z(f)^c \setminus Q}\|f(s)\|^p\,ds $ and $\beta_2=\int\limits_{Q}\|f(s)\|^p\,ds$, then $\beta_1 > \beta_2 >0$.
	Now, define $g\in\;L^p(\mu,X)$ as $g= \beta_1 f \chi_Q - \beta_2 f \chi_{Q^c}$. Then
	\begin{align*}\hspace{-1cm}	\int\limits_{Z(f)^c}\|f(s)\|^{p-1}F_{f(s)}(g(s))\,ds & =\beta_1\int\limits_{Z(f)^c \cap Q}\|f(s)\|^{p} \ ds-\beta_2\int\limits_{Z(f)^c \cap Q^c}\|f(s)\|^{p}\ ds \\
		& = \beta_1 \beta_2 -\beta_2 \beta_1 \\
		&= 0.\end{align*}
	Thus, by \Cref{ourbjp}, $f\bj g$. However, since $Z(f) = Z(g)$, we have
	\begin{align*}
		\hspace{-1.3cm}
		\int\limits_{Z(g)^c}\|g(s)\|^{p-1}F_{g(s)}(f(s))\,ds &=\beta_1^{p-1}\int\limits_{Z(f)^c \cap Q}\|f(s)\|^{p}\,ds-\beta_2^{p-1}\int\limits_{Z(f)^c \cap Q^c}\|f(s)\|^{p}\,ds\\ &=\beta_1^{p-1}\beta_{2}-\beta_2^{p-1}\beta_{1}\neq0,
	\end{align*}
	which gives $g\nbj f$ and hence $f$ is not a left symmetric point.\\
	Next, for $g= \beta_1^{\frac{1}{p-1}} f \chi_Q - \beta_2^{\frac{1}{p-1}} f \chi_{Q^c}\in L^{p}(\mu,X)$, we have
	$$	\int\limits_{Z(g)^c}\|g(s)\|^{p-1}F_{g(s)}(f(s))\,ds =\beta_1\int\limits_{Q}\|f(s)\|^{p}\,ds-\beta_2\int\limits_{Q^c}\|f(s)\|^{p}\,ds=0$$
	and
	\begin{align*}
		\int\limits_{Z(f)^c}\|f(s)\|^{p-1}F_{f(s)}(g(s))\,ds & =\beta_1^{\frac{1}{p-1}}\int\limits_{Q}\|f(s)\|^{p}\,ds-\beta_2^{\frac{1}{p-1}}\int\limits_{Q^c}\|f(s)\|^{p}\,ds\\
		&=\beta_1^{\frac{1}{p-1}}\beta_2-\beta_2^{\frac{1}{p-1}}\beta_1\neq0.
	\end{align*}
	Thus $g\bj f$ and $f\nbj g$ by \Cref{ourbjp}. Hence, $f$ is not a right symmetric point.
\end{proof}
As noticed earlier for the case of $p=1$, the converse of the above result is not true, the same can be checked with \Cref{examatom} and \Cref{examtwoatom}. However, a partial converse of the statement  is true.
\begin{thm}\label{lp}
	Let $X$ be a Banach space whose norm is Fr$\acute{e}$chet differentiable and $f\in L^p(\mu,X),\;1<p<\infty\;p\neq2$ be such that  $Z(f)^c$ is an atom and $f(s)$ is left symmetric (respectively, right symmetric) for each $s\in S$. Then $f$ is left symmetric (respectively, right symmetric).
\end{thm}
\begin{proof}
	Let $g\in L^p(\mu,X)$ be a non-zero element such that $f\bj g$. Then, by \Cref{ourbjp},
	\begin{equation}\label{eq1}
		\int\limits_{Z(f)^c}{\|f(s)\|^{p-1}F_{f(s)}(g(s))\,ds}=	\int\limits_{Z(f)^c\cap Z(g)^c}{\|f(s)\|^{p-1}F_{f(s)}(g(s))\,ds}=0.
	\end{equation}
	Now, $Z(f)^c=A\cup B$, where $A=Z(f)^c\cap Z(g)^c$ and $B=Z(f)^c\cap Z(g)$.  Since $Z(f)^c$ is an atom so either $\mu(A)=0$ or $\mu(B)=0$.
	
	If $\mu(A)=0$ then $$\int\limits_{Z(g)^c}{\|g(s)\|^{p-1}F_{g(s)}(f(s))\,ds}= \int\limits_{A}{\|g(s)\|^{p-1}F_{g(s)}(f(s))\,ds}=0,$$
	so that $g\bj f$ by \Cref{ourbjp}.
	
	If $\mu(B)=0$, write $A=P_1\cup P_2$, where $P_1=\{s\in A: F_{f(s)}(g(s))\geq0\}$ and $P_2=\{s\in A: F_{f(s)}(g(s))<0\}$ both are measurable sets. As done in \Cref{leftsymL1}, if $\mu(P_2)=0$ then Equation (\ref{eq1}) becomes 
	$$\int\limits_{P_1}{\|f(s)\|^{p-1}F_{f(s)}(g(s))\,ds}=0.$$ Thus, $F_{f(s)}(g(s))=0$ for a.e. $s\in P_1$ and hence $F_{g(s)}(f(s))=0$ for a.e. $s\in P_1$,  $f(s)$ being left symmetric point for all $s\in S$. Hence  $$\int\limits_{A}{\|g(s)\|^{p-1}F_{g(s)}(f(s))\,ds}=\int\limits_{P_1}{\|g(s)\|^{p-1}F_{g(s)}(f(s))\,ds}=0.$$ Again by \Cref{ourbjp}, $g\bj f$. The other case when $\mu(P_1)=0$ gives that $\int\limits_{A}{\|f(s)\|^{p-1}F_{f(s)}(g(s))\,ds}=\int\limits_{P_2}{\|f(s)\|^{p-1}F_{f(s)}(g(s))\,ds}\neq0.$ which is a contradiction to $f\bj g$ and hence this case is not possible.  Thus, $f$ is left symmetric.
	
	The right symmetric point can be dealt in the similar manner.
\end{proof}
\begin{cor}
	For a non-atomic measure $\mu$, the space $L^p(\mu, X),\;1\leq p<\infty,\;p\neq2$ has no non-zero left and right symmetric point, where $X$ and $\mu$ are taken same as in \Cref{leftsymL1},  \Cref{rightsymL1} and \Cref{symmLp} respectively.
\end{cor}
For the next consequence, we first provide a  proof of an elementary fact which  was proved for the real Banach spaces in \cite[Proposition 2.1(i)]{paul}.
\begin{lem}\label{RSIL}
	Let $X$ be a Banach space and let $x\in X$ be a non-zero element. If $x$ is right symmteric and smooth then $x$ is left symmteric.
\end{lem}
\begin{proof}
	Let $y\in X$ such that $x\bj y$. By \cite[Theorem 2.1.13]{mal}, there exists $\alpha\in \mathbb{K}$ such that $\alpha x+y\bj x$. Since $x$ is right symmteric, therefore $x\bj \alpha x+y$ and by  \cite[Theorem 2.3.3, Remark 2.3.4]{mal},  $\alpha$ is unique, $x$ being smooth. This gives that $\alpha=0$ as $x\bj y$. Thus, $y\bj x$.
\end{proof}
\begin{cor}\label{symmlp}
	Let $X$ be a Banach space whose norm is Fr$\acute{e}$chet differentiable. If $x=(x_n)\in \ell^{p}(X),\;1 < p<\infty,\; p\neq2$ is a left symmetric point (respectively, right symmetric point) then  either $x=(0,..0,x_{n},0,...)$, where $x_n$ is a left symmetric point (respectively, right symmetric point) in $X$, or $x=(0,..,x_{n},0,x_{m},0,...)$ with $\|x_{n}\|=\|x_{m}\|$. 
\end{cor}
\begin{proof}
	Let $x=(x_n)\in \ell^{p}(X)$ be a left symmetric point, then by \Cref{symmLp}, $Z(x)^c=\{n\in \mathbb{N}: x_n\neq0\}$ is either singleton or $\{n,m\}$.	 If $x=(0,..,0,x_{n},0,...)$,  then as observed in \Cref{l1rsp}, $x_n$ is a left symmetric point in $X$.
	
	Consider the case when $x=(0,..,x_{n},0,x_{m},0,...)$. Let, if possible, $\|x_{n}\|\neq\|x_{m}\|$. For $y=(0,...,0,\frac{\|x_{m}\|^{p-1}}{\|x_{n}\|}x_{n},0,...,0,\frac{-\|x_{n}\|^{p-1}}{\|x_{m}\|}x_{m},0,...)$, we claim that $x\bj y$ but $y\nbj x$. For this, consider 
	\begin{align*}
		\sum_{\{n:\;x_n\neq0\}}{\|x_n\|^{p-1}F_{x_n}(y_n)}&= \|x_{n}\|^{p-1}F_{x_{n}}(y_{n})+\|x_{m}\|^{p-1}F_{x_{m}}(y_{m})\\
		&= \|x_{n}\|^{p-1}\|x_{m}\|^{p-1}-\|x_{n}\|^{p-1}\|x_{m}\|^{p-1}\\
		&= 0
	\end{align*}
	and 
	\begin{align*}
		\sum_{\{n:\;y_n\neq0\}}{\|y_n\|^{p-1}F_{y_n}(x_n)}&= \|y_{n}\|^{p-1}F_{y_{n}}(x_{n})+\|y_{m}\|^{p-1}F_{y_{m}}(x_{m})\\
		&= \|y_{n}\|^{p-1}\|x_{n}\|-\|y_{m}\|^{p-1}\|x_{m}\|\\
		&= \|x_{m}\|^{(p-1)^2}\|x_{n}\|-\|x_{n}\|^{(p-1)^2}\|x_{m}\|\\
		&\neq0,
	\end{align*} 
	by using the fact that $p\neq2$ and $\|x_{n}\|\neq\|x_{m}\|$. Hence, by \Cref{ourbjp}, $x\bj y$ and $y\nbj x$ which is a contradiction. This proves the claim.
	
	Note that $\ell^{p}(X)$ is smooth, $X$ being smooth. Therefore, by \Cref{RSIL}, every non-zero right symmetric point is also left symmetric point. Further, it is easy to check that if $(0,..,0,x_n,0,...)$ is right symmetric then $x_n$ is right symmetric in $X$ and  thus the result.
\end{proof}
Converse of the above result need not be true. In fact, we prove that an element of the latter form cannot be left symmteric (respectively, right symmetric) if $X$ is a Hilbert space of dimension more than 1. 
\begin{cor}\label{lpcharc}
	For a  Hilbert space $H$ with dim$(H)\geq2$, an element $x\in \ell^{p}(H),\;1 < p<\infty,\; p\neq2$, is left symmetric (respectively, right symmetric) if and only if $x=(0,..0,h,0,...)$ for some $h \in H$.
\end{cor}
\begin{proof}
	In view of \Cref{symmlp} and \Cref{lp}, it is sufficient to show that  
	an element $x$ of the form $(0,..,x_{n},0,x_{m},0,...)$ with $0 \neq \|x_{n}\|=\|x_{m}\|$ is not a left symmetric point in $\ell^{p}(H)$.
	Since $\dim(H)\geq2$, there exists a non-zero $x_{0}\in H$ such that $\langle x_{n},x_0\rangle=0$. Consider $y=(0,....0, y_n,0,...0,y_{m},0,....) \in \ell^{p}(H)$, where $y_n=-x_{n}+x_0,$ and $y_m=x_m$. Since $\|x_{n}\|=\|x_{m}\|$, we have $\langle y_{n},x_{n}\rangle+\langle y_{m},x_{m}\rangle=0$ and
	\begin{align*}
		\|x_{n}\|^{p-1}F_{x_{n}}(y_{n})+\|x_{m}\|^{p-1}F_{x_{m}}(y_{m}) &=\|x_{n}\|^{p-2}\langle y_{n},x_{n}\rangle+\|x_{m}\|^{p-2}\langle y_{m},x_{m}\rangle \\
		& = \|x_{n}\|^{p-2}(\langle y_{n},x_{n}\rangle+ \langle y_{m},x_{m}\rangle)=0.
	\end{align*}
	Thus, by \Cref{ourbjp}, $x\bj y$. However,  
	\begin{align*}
		\|y_{n}\|^{p-1}F_{y_{n}}(x_{n})+\|y_{m}\|^{p-1}F_{y_{m}}(x_{m}) & = \|y_{n}\|^{p-2}\langle x_{n},y_{n}\rangle+\|y_{m}\|^{p-2}\langle x_{m},y_{m}\rangle \\
		& = (-1)(\|y_{n}\|^{p-2}-\|y_{m}\|^{p-2})
	\end{align*}
	Again, by  \Cref{ourbjp},  $y\bj x$ if and only if $\|y_{n}\|^{p-2}=\|y_{m}\|^{p-2}$, or equivalently, $x_0 = 0$, which is not true. Thus, $y\nbj x$ and hence $x$ is not a left symmetric point.
\end{proof}	 
 \begin{remark}
	In the last result it is evident that if $f$ is left (respectively, right) symmetric then $f(s)$ is left (respectively, right) symmetric  for every $s \in S$ . It would be interesting to know whether  this holds true in $L^p(\mu, X)$.  
\end{remark}

\end{document}